\documentclass{article}
\usepackage{amssymb}
\usepackage{amsfonts}
\usepackage{amsmath}

\begin{document}

\begin{center}
{\Large Some applications of Fibonacci and Lucas numbers }

\begin{equation*}
\end{equation*}%
Cristina FLAUT, Diana SAVIN and Gianina ZAHARIA 
\begin{equation*}
\end{equation*}
\end{center}

\textbf{Abstract. }{\small In this paper, we provide new applications of
Fibonacci and Lucas numbers. In some circumstances, we find algebraic
structures on some sets defined with these numbers, we generalize Fibonacci
and Lucas numbers by using an arbitrary binary relation over the real fields
instead of addition of the real numbers and we give properties of the new
obtained sequences. Moreover, by using some relations between Fibonacci and
Lucas numbers, we provide a method to find new examples of split quaternion
algebras and we give new properties of these elements. }%
\begin{equation*}
\end{equation*}


\medskip

\begin{equation*}
\end{equation*}

\bigskip \textbf{1. Introduction}%
\begin{equation*}
\end{equation*}

There have been numerous papers devoted to the study of the properties and
applications of Fibonacci and Lucas sequences. Due to this fact, to obtain
new results in this direction is not always an easy problem. Our purpose in
this chapter is to provide new applications of these sequences. Therefore,
we studied if some algebraic structures can be defined by using Fibonacci
and Lucas elements and we provided some new properties and relations of them.

Moreover, by taking into consideration some known relations between these
numbers, we give a method to find new examples of quaternion split algebras.
This method can be an easy alternative to using the properties of quadratic
forms.

Let $\left( f_{n}\right) _{n\geq 0}$ be the Fibonacci sequence 
\begin{equation*}
f_{n}=f_{n-1}+f_{n-2},\;n\geq 2,f_{0}=0;f_{1}=1,
\end{equation*}%
and $\left( l_{n}\right) _{n\geq 0}$ be the Lucas sequence 
\begin{equation*}
l_{n}=l_{n-1}+l_{n-2},\;n\geq 2,l_{0}=2;l_{1}=1.
\end{equation*}%
Let $n$ be an arbitrary positive integer and let $p,q$ be two arbitrary
integers. In the paper [Fl, Sa; 15], we introduced the sequence $\left(
g_{n}\right) _{n\geq 1},$ where 
\begin{equation}
g_{n+1}=pf_{n}+ql_{n+1},\;n\geq 0,  \tag{1.1}
\end{equation}%
with $g_{0}=p+2q,g_{1}=q$. We remark that $g_{n}=g_{n-1}+g_{n-2}.$\newline
This sequence is called \textit{\ generalized Fibonacci-Lucas numbers}. To
avoid confusions, from now on, we will use the notation $g_{n}^{p,q}$
instead of $g_{n}.$\newline

Let $\mathbb{H}\left( \alpha ,\beta \right) $ be the generalized real\
quaternion algebra, i.e. the algebra of the elements of the form $%
a=a_{1}\cdot 1+a_{2}e_{2}+a_{3}e_{3}+a_{4}e_{4},$ where $a_{i}\in \mathbb{R}%
,i\in \{1,2,3,4\}$, and the elements of the basis $\{1,e_{2},e_{3},e_{4}\}$
satisfy the following rules, given in the below multiplication
table:\medskip\ \vspace{3mm}

\begin{center}
\begin{tabular}{c|cccc}
$\cdot $ & $1$ & $e_{1}$ & $e_{2}$ & $e_{3}$ \\ \hline
$1$ & $1$ & $e_{1}$ & $e_{2}$ & $e_{3}$ \\ 
$e_{1}$ & $e_{1}$ & $\alpha $ & $e_{3}$ & $\alpha e_{2}$ \\ 
$e_{2}$ & $e_{2}$ & $-e_{3}$ & $\beta $ & $-\beta e_{1}$ \\ 
$e_{3}$ & $e_{3}$ & $-\alpha e_{2}$ & $\beta e_{1}$ & $-\alpha \beta $%
\end{tabular}%
.\medskip 
\end{center}

We denote by $\boldsymbol{t}\left( a\right) $ and $\mathbf{n}\left( a\right) 
$ \textit{the trace} and \textit{the norm} of a quaternion $a$. The norm of
a generalized quaternion has the following expression 
\begin{equation}
\mathbf{n}\left( a\right) =a_{1}^{2}-\alpha a_{2}^{2}-\beta a_{3}^{2}+\alpha
\beta a_{4}^{2}.  \tag{1.2.}
\end{equation}%
If, \thinspace for $x\in \mathbb{H}\left( \alpha ,\beta \right) $, the
relation $\mathbf{n}\left( x\right) =0$ implies $x=0$, then the algebra $%
\mathbb{H}\left( \alpha ,\beta \right) $ is called a \textit{division}
algebra, otherwise the quaternion algebra is called a \textit{split}
algebra. For $\alpha =\beta =-1$, we obtain the real division algebra $%
\mathbb{H}$.

\begin{equation*}
\end{equation*}

\bigskip \textbf{2. Preliminaries}

\begin{equation*}
\end{equation*}

First of all, we recall some elementary properties of the Fibonacci and
Lucas numbers, properties which will be used in this chapter.\newline
Let $(f_{n})_{n\geq 0}$ be the Fibonacci sequence and let $(l_{n})_{n\geq 0}$
be the Lucas sequence. Let $\alpha =\frac{1+\sqrt{5}}{2}$ and $\beta =\frac{%
1-\sqrt{5}}{2}$, be two real numbers.

The following formulae are well known:\medskip \newline
\smallskip \textit{Binet's formula for Fibonacci sequence} 
\begin{equation*}
f_{n}=\frac{\alpha ^{n}-\beta ^{n}}{\alpha -\beta }=\frac{\alpha ^{n}-\beta
^{n}}{\sqrt{5}},\ \ n\in \mathbb{N}.
\end{equation*}%
\textit{Binet's formula for Lucas sequence} 
\begin{equation*}
l_{n}=\alpha ^{n}+\beta ^{n},\ \ n\in \mathbb{N}.
\end{equation*}

\medskip

\textbf{Proposition 2.1.} ([Fib.]). \textit{Let} $(f_{n})_{n\geq 0}$ \textit{%
be the Fibonacci sequence} \textit{and let } $(l_{n})_{n\geq 0}$ \textit{be
the Lucas sequence.} \textit{Therefore, the following properties hold:}%
\newline
i) 
\begin{equation*}
5\mid f_{n}\;\text{\textit{if\ and\ only\ if}}\;5\mid n.
\end{equation*}%
ii) 
\begin{equation*}
f_{2n}=l_{n}\cdot f_{n},\ \ n\in \mathbb{N}.
\end{equation*}%
iii) 
\begin{equation*}
f_{n}^{2}+f_{n+1}^{2}=f_{2n+1},\ \ n\in \mathbb{N}.
\end{equation*}%
\newline
iv)%
\begin{equation*}
f_{n}^{2}-f_{n+1}f_{n-1}=\left( -1\right) ^{n-1},n\in \mathbb{N}.
\end{equation*}%
\newline
v)%
\begin{equation*}
f_{-n}=\left( -1\right) ^{n+1}f_{n},n\in \mathbb{N}.
\end{equation*}%
\newline
vi) 
\begin{equation*}
l_{-n}=\left( -1\right) ^{n}l_{n},n\in \mathbb{N}.
\end{equation*}%
\newline
vii) 
\begin{equation*}
l_{n}^{2}=5f_{n}^{2}+4\left( -1\right) ^{n},n\in \mathbb{N}.
\end{equation*}%
\newline
viii)%
\begin{equation*}
l_{2n}=l_{n}^{2}+2\left( -1\right) ^{n+1},n\in \mathbb{N}..
\end{equation*}%
\newline
ix) 
\begin{equation*}
l_{2n}l_{2n+2}-5f_{n+1}^{2}=1,n\in \mathbb{N}.
\end{equation*}%
\newline
x)\newline
\begin{equation*}
f_{2n}+f_{n}^{2}=2f_{n}f_{n+1},n\in \mathbb{N}.
\end{equation*}%
\newline
xi)%
\begin{equation*}
f_{2n}-f_{n}^{2}=2f_{n}f_{n-1},n\in \mathbb{N}.
\end{equation*}%
\newline
xii)%
\begin{equation*}
l_{n}^{2}-f_{n}^{2}=4f_{n-1}f_{n+1},n\in \mathbb{N}.
\end{equation*}%
\newline
xiii)\smallskip \newline
\begin{equation*}
f_{2n}=f_{n+1}^{2}-f_{n-1}^{2},n\in \mathbb{N}.
\end{equation*}%
xiv)\smallskip 
\begin{equation*}
l_{4n}=5f_{2n}^{2}+2,n\in \mathbb{N}.
\end{equation*}%
$\Box \medskip $

\textbf{Proposition 2.2.} ([Gi, Sz; 06]). \textit{Let} $K$ \textit{be a
field. Then, the quaternion algebra} $\mathbb{H}_{K}\left( a,b\right) $ 
\textit{is a split algebra if and only if the conic} 
\begin{equation*}
C\left( a,b\right) :ax^{2}+by^{2}=z^{2}
\end{equation*}%
\textit{has a rational point over} $K^{\ast }$, \textit{i.e. there are} $%
x_{0},y_{0},z_{0}\in K^{\ast }$ \textit{such that} $%
ax_{0}^{2}+by_{0}^{2}=z_{0}^{2}$. $\Box $

\begin{equation*}
\end{equation*}%
\textbf{3. Main results}

\begin{equation*}
\end{equation*}

For the beginning, we want to find some properties of Fibonacci numbers $%
f_{5n},$ $n\in \mathbb{N}.$ We start with some examples. \medskip

\textbf{Example 3.1.} $f_{10}=55=11f_{5},$ $%
f_{15}=610=11f_{10}+f_{5}=122f_{5},$ $f_{20}=11f_{15}+f_{10}=6765=1353f_{5}.$%
\medskip

Using the computer algebra system Magma ([Mag.]), by using the function 
\textit{Fibonacci(n)}, we obtain another example.

We get $f_{65}=17167680177565,$ $f_{35}=9227465$ and $f_{65}$ \textit{div} $%
f_{35}$ is $1860497,$ $f_{65}$ \textit{mod} $f_{35}$ is $9227460$ and $%
9227460$ \textit{div} $f_{5}$ is $1845492.$ Therefore, it results that 
\begin{equation*}
f_{65}=1860497f_{35}+9227460=1860497f_{35}+1845492f_{5}.
\end{equation*}

\textbf{Proposition 3.2.} \textit{The set} 
\begin{equation*}
A=\left\{ \alpha f_{5n}|n\in \mathbb{N},\alpha \in \mathbb{Z}\right\}
\end{equation*}%
\textit{is a commutative non-unitary ring, with addition and
multiplication.\medskip }

\textbf{Proof.} We remark that $f_{0}=0\in A$ is the identity element for
addition on $A$. It is clear that the addition on $A$ is commutative.\newline
We consider the Fibonacci number $f_{5n}.$ If $n\neq 0$, there are two cases.%
\newline
\textbf{Case 1.} If $n$ is even, then $n=2k$, where $k\in \mathbb{N}-\{0\}$.
Applying Proposition 2.1 ii), we have that $f_{5n}=f_{10k}=l_{5k}f_{5k}.$ If
we denote $l_{5k}=a\in \mathbb{N}-\{0\}$, we obtain $f_{5n}=f_{10k}=af_{5k}$
and, applying Proposition 2.1 i), we obtain $f_{5k}=a^{\prime }f_{5}$, where 
$a^{^{\prime }}\in \mathbb{N}-\{0\}$. It results, 
\begin{equation}
f_{5n}=f_{10k}=af_{5k}=aa^{^{\prime }}f_{5},\ \text{with}\ a,a^{^{\prime
}}\in \mathbb{N}-\{0\}.  \tag{3.1}
\end{equation}%
\textbf{Case 2. }If $n$ is odd, then $n=2k+1$, where $k\in \mathbb{N}-\{0\}$%
. We have 
\begin{equation*}
f_{5n}=f_{10k+5}=f_{10k+4}+f_{10k+3}=2f_{10k+3}+f_{10k+2}=
\end{equation*}%
\begin{equation*}
=3f_{10k+2}+2f_{10k+1}=5f_{10k+1}+3f_{10k}=f_{10k+1}f_{5}+3f_{10k}.
\end{equation*}%
If we denote $f_{10k+1}=a$$\in \mathbb{N}-\{0\},$ we obtain $%
f_{5n}=af_{5}+3f_{10k}.$ Applying Proposition 2.1 i), we get $f_{10k}=bf_{5},
$ where $b$$\in \mathbb{N}-\{0\}$. It results 
\begin{equation}
f_{5n}=f_{10k+5}=af_{5}+3f_{10k}=\left( a+3b\right) f_{5},\ \text{with}\
a,b\in \mathbb{N}-\{0\}.  \tag{3.2}
\end{equation}%
From the relations (3.1) and (3.2), we have that $\alpha f_{5n}-\beta
f_{5m}=\gamma f_{5l}$, for each $n,m\in \mathbb{N}$ and for each $\alpha
,\beta \in \mathbb{Z}$,where $l\in \mathbb{N},$ $\gamma \in \mathbb{Z}$.
Therefore $\left( A,+\right) $ is a subgroup of the group $\left( \mathbb{Z}%
,+\right) .\newline
$Applying Proposition 2.1 i), we obtain that $\alpha f_{5n}\beta
f_{5m}=\delta f_{5r}$, for each $n,m\in \mathbb{N}$ and for each $\alpha
,\beta \in \mathbb{Z}$, where $r\in \mathbb{N},$ $\delta \in \mathbb{Z}$.
From here, it is clear that the multiplication on $A$ is commutative.\newline
Therefore $\left( A,+,\cdot \right) $ is a commutative non-unitary subring
of the ring $\left( \mathbb{Z},+,\cdot \right) $. $\Box \medskip $

Using the above proposition, we obtain the below results.\medskip

\textbf{Proposition 3.3.} \textit{We consider} $\left( g_{n}^{p,q}\right)
_{n\geq 1}$\textit{the generalized Fibonacci-Lucas numbers and the set }%
\begin{equation*}
M=\left\{ \sum\limits_{i=1}^{n}g_{5n_{i}}^{p_{i},5q_{i}}|n_{i}\in \mathbb{N}%
,p_{i},q_{i}\in \mathbb{Z}\right\} \cup \left\{ 0\right\} .
\end{equation*}%
\textit{Then, the following relations are true:\newline
}i) $\left( M,+\right) $ \textit{is a subgroup of the group} $\left( 5%
\mathbb{Z},+\right) .$\newline
ii) $M$ \textit{is a bilateral ideal of the ring} $\left( \mathbb{Z},+,\cdot
\right) $.\medskip \medskip

\textbf{Proof.} Using Proposition 2.1 i), we have $5\mid g_{5n}^{\alpha
p,5\alpha q}$, for all$~n\in \mathbb{N}$ and $p,q\in \mathbb{Z}$. Let $%
n,m\in \mathbb{N}-\{0\}$ and $p,q,p^{^{\prime }},q^{^{\prime }},\alpha
,\beta \in \mathbb{Z}$. We obtain that 
\begin{equation*}
\alpha g_{5n}^{p,5q}+\beta g_{5m}^{p^{^{\prime }},5q^{^{\prime
}}}=g_{5n}^{\alpha p,5\alpha q}+g_{5m}^{\alpha p^{^{\prime }},5\alpha
q^{^{\prime }}}.
\end{equation*}%
From here, we obtain conditions i) and ii). $\Box \medskip $\medskip \newline

\textbf{Proposition 3.4.} \textit{With the above notations, the following
statements are true:}\newline
1) \textit{The quaternion algebra} $\mathbb{H}_{\mathbb{Q}}\left(
f_{10n+5},-1\right) $ \textit{is a split algebra over} $\mathbb{Q}$, $n\in 
\mathbb{N}.$\newline
2) \textit{The quaternion algebra} $\mathbb{H}_{\mathbb{Q}}\left( \frac{1}{5}%
l_{20n},-\frac{2}{5}\right) $ \textit{is a split algebra over} $\mathbb{Q}$, 
$n\in \mathbb{N}.$\newline
3) \textit{The quaternion algebra} $\mathbb{H}_{\mathbb{Q}}\left(
f_{n+1}f_{n-1},\left( -1\right) ^{n-1}\right) $ \textit{is a split algebra
over} $\mathbb{Q}$, $n\in \mathbb{N}.$ \newline
4) \textit{The quaternion algebra }$\mathbb{H}_{\mathbb{Q}}\left( 5,\left(
-1\right) ^{n}\right) $ \textit{is a split algebra over} $\mathbb{Q}$, $n\in 
\mathbb{N}.$\newline
5) \textit{The quaternion algebra }$\mathbb{H}_{\mathbb{Q}}\left(
l_{2n}l_{2n+2},-5\right) $ \textit{is a split algebra over} $\mathbb{Q}$, $%
n\in \mathbb{N}.$\newline
6) \textit{The quaternion algebra }$\mathbb{H}_{\mathbb{Q}}\left(
2f_{n}f_{n+1},-f_{2n}\right) $ \textit{is a split algebra over} $\mathbb{Q}$%
, $n\in \mathbb{N}.$\newline
7) \textit{The quaternion algebra }$\mathbb{H}_{\mathbb{Q}}\left(
f_{2n},-2f_{n}f_{n-1}\right) $ \textit{is a split algebra over} $\mathbb{Q}$%
, $n\in \mathbb{N}.$\newline
8) \textit{The quaternion algebras }$\mathbb{H}_{\mathbb{Q}}\left(
f_{n-1}f_{n+1},f_{n}^{2}\right) $ \textit{and} $\mathbb{H}_{\mathbb{Q}%
}\left( 1,-f_{n-1}f_{n+1}\right) $ \textit{are split algebras over} $\mathbb{%
Q}$, $n\in \mathbb{N}.$\newline
9) \textit{The quaternion algebra }$\mathbb{H}_{\mathbb{Q}}\left(
f_{2n},1\right) $ \textit{is a split algebra over} $\mathbb{Q}$, $n\in 
\mathbb{N}$.\medskip 

\textbf{Proof.} 1) Applying Proposition 3.4 (i) from [Sa; 19], for $%
n\rightarrow 5n+2,$ we obtain that the quaternion algebra $\mathbb{H}_{%
\mathbb{Q}}\left( f_{10n+5},-1\right) $ is a split algebra over $\mathbb{Q}$%
, for all $n\in \mathbb{N}$.\newline
2) For this purpose, we use Proposition 2.2. Therefore, we study if the
equation 
\begin{equation}
\frac{1}{5}l_{20n}\cdot x^{2}-\frac{2}{5}\cdot y^{2}=z^{2}  \tag{3.3.}
\end{equation}%
has rational solutions. From Proposition 2.1 xiv), we have $%
l_{20n}=5f_{10n}^{2}+2$, for all$~n\in \mathbb{N}$. It results $\frac{1}{5}%
l_{20n}-\frac{2}{5}=f_{10n}^{2}$, then $\left( x_{0},y_{0},z_{0}\right)
=\left( 1,1,f_{10n}\right) $$\in $$\mathbb{Q}\times \mathbb{Q}\times \mathbb{%
Q}$, is a rational solution of equation $\left( 3.3\right) $. We obtain that
the quaternion algebra $\mathbb{H}_{\mathbb{Q}}\left( \frac{1}{5}l_{20n},-%
\frac{2}{5}\right) $ is a split algebra over $\mathbb{Q}$, for all $n\in 
\mathbb{N}$.\newline
3) Using the same idea, we study if the equation 
\begin{equation}
f_{n+1}f_{n-1}x^{2}+\left( -1\right) ^{n-1}y^{2}=z^{2}  \tag{3.4.}
\end{equation}%
has rational solutions. From Proposition 2.1 iv), it results that $%
f_{n}^{2}=f_{n+1}f_{n-1}+\left( -1\right) ^{n-1}$, for all $n\in \mathbb{N}$%
, therefore $\left( 1,1,f_{n}\right) $ is a rational solution of equation $%
\left( 3.4\right) $. We obtain that the quaternion algebra $\mathbb{H}_{%
\mathbb{Q}}\left( f_{n+1}f_{n-1},\left( -1\right) ^{n-1}\right) $ is a split
algebra over $\mathbb{Q}$, for all $n\in \mathbb{N}$.\newline
4) From Proposition 2.1 vii), we have $l_{n}^{2}=5f_{n}^{2}+4\left(
-1\right) ^{n}$and it results that the equation $5x^{2}+\left( -1\right)
^{n-1}y^{2}=z^{2}$ has $\left( f_{n},2,l_{n}\right) $ as a rational
solution. Therefore, the quaternion algebra $\mathbb{H}_{\mathbb{Q}}\left(
5,\left( -1\right) ^{n}\right) $ is a split algebra over $\mathbb{Q}$, for
all $n\in \mathbb{N}$.\newline
5) We study if the equation 
\begin{equation}
l_{2n}l_{2n+2}x^{2}-5y^{2}=z^{2}  \tag{3.5.}
\end{equation}%
has rational solutions. From Proposition 2.1 ix), we have $%
l_{2n}l_{2n+2}-5f_{n+1}^{2}=1,~n\in \mathbb{N}$, it results that the pair $%
\left( 1,f_{n+1}^{2},1\right) $ is a rational solution of the equation $%
\left( 3.5\right) $. Therefore, the quaternion algebra$~\mathbb{H}_{\mathbb{Q%
}}\left( l_{2n}l_{2n+2},-5\right) $ is a split algebra over $\mathbb{Q}$,
for all $n\in \mathbb{N}$.\newline
6) We study if the equation 
\begin{equation}
2f_{n}f_{n+1}x^{2}-f_{2n}y^{2}=z^{2}  \tag{3.6.}
\end{equation}%
has rational solutions. From Proposition 2.1 x), we have that the pair $%
\left( 1,1,f_{n}\right) $ is a rational solution of the equation $\left(
3.6\right) $. Therefore, the quaternion algebra $\mathbb{H}_{\mathbb{Q}%
}\left( 2f_{n}f_{n+1},-f_{2n}\right) $ is a split algebra over $\mathbb{Q}$,
for all $n\in \mathbb{N}$.\newline
7) From Proposition 2.1, xi), the equation 
\begin{equation}
f_{2n}x^{2}-2f_{n}f_{n-1}y^{2}=z^{2}  \tag{3.7.}
\end{equation}%
\newline
has the rational solution $\left( 1,1,f_{n}\right) $, therefore the
quaternion algebra $\mathbb{H}_{\mathbb{Q}}\left(
f_{2n},-2f_{n}f_{n-1}\right) $ is a split algebra over $\mathbb{Q}$, for all 
$n\in \mathbb{N}$.\newline
8) From Proposition 2.1 xii), the equation%
\begin{equation}
f_{n-1}f_{n+1}x^{2}+f_{n}^{2}y^{2}=z^{2}  \tag{3.8.}
\end{equation}%
has the rational solution $\left( 2,1,l_{n}\right) $, therefore the
quaternion algebra $\mathbb{H}_{\mathbb{Q}}\left(
f_{n-1}f_{n+1},f_{n}^{2}\right) $ is a split algebra over $\mathbb{Q}$, for
all $n\in \mathbb{N}$. \newline
In the same way, the equation 
\begin{equation}
x^{2}-f_{n-1}f_{n+1}y^{2}=z^{2}  \tag{3.9.}
\end{equation}%
has the rational solution $\left( l_{n},2,f_{n}\right) $, therefore the
quaternion algebra $\mathbb{H}_{\mathbb{Q}}\left( 1,-f_{n-1}f_{n+1}\right) $
is a split algebra over $\mathbb{Q}$, for all $n\in \mathbb{N}$. \newline
9) From Proposition 2.1 xiii), the equation%
\begin{equation}
f_{2n}x^{2}+y^{2}=f_{n+1}^{2}  \tag{3.10.}
\end{equation}%
has the rational solution $\left( 1,f_{n-1},f_{n+1}\right) $, therefore the
quaternion algebra $\mathbb{H}_{\mathbb{Q}}\left( f_{2n},1\right) $ is a
split algebra over $\mathbb{Q}$, for all $n\in \mathbb{N}$. $\Box \medskip $%
\medskip \newline

Let $\mathbb{H}\left( \alpha ,\beta \right) $ be the generalized real\
quaternion algebra with the basis $\{1,e_{2},e_{3},e_{4}\}$, $(f_{n})_{n\geq
0}$ be the Fibonacci sequence and $(l_{n})_{n\geq 0}$ be the Lucas sequence.
We consider the quaternions 
\begin{equation*}
F_{n}=f_{n}+f_{n+1}e_{2}+f_{n+3}e_{3}+f_{n+4}e_{4}
\end{equation*}%
and%
\begin{equation*}
L_{n}=l_{n}+l_{n+1}e_{2}+l_{n+3}e_{3}+l_{n+4}e_{4},
\end{equation*}%
called \textit{Fibonacci quaternion} and, respectively, \textit{Lucas
quaternion}.\medskip 

\textbf{Proposition 3.5.} \textit{With the above notations, in }$\mathbb{H}%
\left( \alpha ,\beta \right) $\textit{, the following relations are true:}

1) $5\mathbf{n}\left( F_{n}\right) =\mathbf{n}\left( L_{n}\right) ,n\in 
\mathbb{N}.$

2) $5(f_{2n+1}+f_{2n+5})=l_{2n}+l_{2n+2}+l_{2n+4}+l_{2n+6},n\in \mathbb{N}.$

3) $\mathbf{n}\left( F_{n}+L_{n}\right) =\mathbf{n}\left( F_{n})+\mathbf{n}%
(L_{n}\right) +2\left( f_{2n+7}-f_{2n-1}\right) ,n\in \mathbb{N}$.\medskip 

\textbf{Proof.} 1) From Proposition 2.1 vii), we have\newline
$\mathbf{n}(L_{n})=l_{n}^{2}+l_{n+1}^{2}+l_{n+3}^{2}+l_{n+4}^{2}=$\newline
$=5f_{n}^{2}+4\left( -1\right) ^{n}+5f_{n+1}^{2}+4\left( -1\right) ^{n+1}+$%
\newline
$+5f_{n+2}^{2}+4\left( -1\right) ^{n+2}+5f_{n+3}^{2}+4\left( -1\right)
^{n+3}=5\mathbf{n}\left( F_{n}\right) $.

2) From Proposition 2.1 viii), we have\newline
$\mathbf{n}(L_{n})=l_{n}^{2}+l_{n+1}^{2}+l_{n+3}^{2}+l_{n+4}^{2}=$\newline
$=l_{2n}-$ $2\left( -1\right) ^{n+1}+l_{2n+2}-$ $2\left( -1\right) ^{n+2}+$%
\newline
$+l_{2n+4}-$ $2\left( -1\right) ^{n+3}+l_{2n+6}-$ $2\left( -1\right) ^{n+4}=$%
\newline
$=l_{2n}+l_{2n+2}+l_{2n+4}+l_{2n+6}.$

From Proposition 2.1 iii),we have\newline
$\mathbf{n}%
(F_{n})=f_{n}^{2}+f_{n+1}^{2}+f_{n+2}^{2}+f_{n+3}^{2}=f_{2n+1}+f_{2n+5}$ and
we apply the relation 1).

3) We have $\mathbf{n}\left( F_{n}+L_{n}\right) =$\newline
$%
=(f_{n}+l_{n})^{2}+(f_{n+1}+l_{n+1})^{2}+(f_{n+2}+l_{n+2})^{2}+(f_{n+3}+l_{n+3})^{2}= 
$\newline
$%
=f_{n}^{2}+f_{n+1}^{2}+f_{n+2}^{2}+f_{n+3}^{2}+l_{n}^{2}+l_{n+1}^{2}+l_{n+3}^{2}+l_{n+4}^{2}+ 
$\newline
$+2f_{n}l_{n}+2f_{n+1}l_{n+1}+2f_{n+2}l_{n+2}+2f_{n+3}l_{n+3}=$\newline
$=\mathbf{n}\left( F_{n})+\mathbf{n}(L_{n}\right) +2\left(
f_{2n}+f_{2n+2}+f_{2n+4}+f_{2n+6}\right) =$\newline
$=\mathbf{n}\left( F_{n})+\mathbf{n}(L_{n}\right) +2\left(
f_{n+1}^{2}-f_{n-1}^{2}+f_{n+2}^{2}-f_{n}^{2}+f_{n+3}^{2}-f_{n+1}^{2}+f_{n+4}^{2}-f_{n+2}^{2}\right) = 
$\newline
$=\mathbf{n}\left( F_{n})+\mathbf{n}(L_{n}\right) +2\left(
-f_{2n-1}+f_{2n+7}\right) $. $\Box \medskip $\medskip \newline

Fibonacci numbers are defined by using the addition operation over the real
field. Will be interesting to search what is happen when, instead of
addition, we use an arbitrary binary relation over $\ \mathbb{R}$.
Therefore, we will consider \textquotedblright $\ast $\textquotedblright\ a
binary relation over $\mathbb{R}$ and $a,b\in \mathbb{R}.$ We define the
following sequence 
\begin{equation*}
\varphi _{n}=\varphi _{n-1}\ast \varphi _{n-2},\varphi _{0}=a,\varphi
_{1}=b.~\ 
\end{equation*}%
We call this sequence the \textit{left} \textit{Fibonacci type sequence
attached to the binary relation }\textquotedblright $\ast $%
\textquotedblright ,~\textit{generated by }$a$ \textit{and} $b$.

The following sequence 
\begin{equation*}
\varphi _{n}=\varphi _{n-2}\ast \varphi _{n-1},\varphi _{0}=a,\varphi
_{1}=b~\ 
\end{equation*}%
is called the \textit{right} \textit{Fibonacci type sequence attached to the
binary relation }\textquotedblright $\ast $\textquotedblright , \textit{\
generated by }$a$ \textit{and} $b$.\medskip

\textbf{Proposition 3.6.} \textit{Let} $A,B\in \mathbb{R}$ \textit{such that}
$\Delta =A^{2}+4B>0$. \textit{We define on} $\mathbb{R}$ \textit{the
following binary relation} 
\begin{equation*}
x\ast y=Ax+By,x,y\in \mathbb{R}\text{.}
\end{equation*}%
\textit{We consider} $\left( \varphi _{n}\right) _{n\in N}$ \textit{the left
Fibonacci type sequence attached to the binary relation \textquotedblright }$%
\ast $\textquotedblright \textit{. Therefore, we have}%
\begin{equation}
\varphi _{n+1}=\frac{1}{\beta -\alpha }\left[ \left( -b+a\beta \right)
\alpha ^{n+1}+\left( b-a\alpha \right) \beta ^{n+1}\right] ,  \tag{3.11.}
\end{equation}%
\textit{where} $\alpha =\frac{A+\sqrt{\Delta }}{2},\beta =\frac{A-\sqrt{%
\Delta }}{2}$ \textit{and} 
\begin{equation}
\mathcal{L}=\lim \frac{\varphi _{n+1}}{\varphi _{n}}=\max \{\alpha ,\beta \}.
\tag{3.12.}
\end{equation}

\textbf{Proof.} Taking $\varphi _{0}=a,\varphi _{1}=b$, we get\newline
$\varphi _{2}=\varphi _{1}\ast \varphi _{0}=Ab+Ba,$\newline
$\varphi _{3}=\varphi _{2}\ast \varphi _{1}=A\left( Ab+Ba\right) +Bb=$%
\newline
$=\left( A^{2}+B\right) b+ABa$, etc. From here, $\varphi _{n}$ can be write
under the form 
\begin{equation*}
\varphi _{n}=\varphi _{n-1}\ast \varphi _{n-2}=x_{n}b+y_{n}a.
\end{equation*}%
We have%
\begin{equation*}
\varphi _{n+1}=x_{n+1}b+y_{n+1}a
\end{equation*}%
and\newline
$\varphi _{n+1}=A\varphi _{n}+B\varphi _{n-1}=A\left( x_{n}b+y_{n}a\right)
+B\left( x_{n-1}b+y_{n-1}a\right) =$\newline
$=\left( Ax_{n}+Bx_{n-1}\right) b+\left( A_{y_{n}}+By_{n-1}\right) a.$ It
results, 
\begin{equation*}
x_{n+1}=Ax_{n}+Bx_{n-1}
\end{equation*}%
and%
\begin{equation*}
y_{n+1}=Ay_{n}+By_{n-1.}
\end{equation*}%
We consider the roots of the equation 
\begin{equation*}
x^{2}-Ax-B=0,
\end{equation*}%
\begin{equation*}
\alpha =\frac{A+\sqrt{\Delta }}{2},\beta =\frac{A-\sqrt{\Delta }}{2},\Delta
=A^{2}+4B.
\end{equation*}%
Therefore, 
\begin{equation*}
x_{n+1}=X_{1}\alpha ^{n+1}+X_{2}\beta ^{n+1},
\end{equation*}%
\begin{equation*}
y_{n+1}=Y_{1}\alpha ^{n+1}+Y_{2}\beta ^{n+1},
\end{equation*}%
with $X_{1}+X_{2}=0,X_{1}\alpha +X_{2}\beta =1$ and $Y_{1}+Y_{2}=1,Y_{1}%
\alpha +Y_{2}\beta =0$. We obtain $X_{1}=-\frac{1}{\beta -a},X_{2}=\frac{1}{%
\beta -a},Y_{1}=\frac{\beta }{\beta -a},Y_{2}=-\frac{\alpha }{\beta -a}$. It
results, 
\begin{equation*}
x_{n+1}=\frac{1}{\beta -a}(-\alpha ^{n+1}+\beta ^{n+1})
\end{equation*}%
and 
\begin{equation*}
y_{n+1}=\frac{1}{\beta -a}(\beta \alpha ^{n+1}+\alpha \beta ^{n+1}).
\end{equation*}%
Therefore, by straightforward calculations, we get 
\begin{equation*}
\varphi _{n+1}=\frac{1}{\beta -\alpha }\left[ \left( -b+a\beta \right)
\alpha ^{n+1}+\left( b-a\alpha \right) \beta ^{n+1}\right] .
\end{equation*}%
This formula is a Binet-type formula.

Now, for the limit $\mathcal{L}$, we obtain%
\begin{equation*}
\mathcal{L}\text{=}\underset{n\rightarrow \infty }{\lim }\frac{\varphi _{n+1}%
}{\varphi _{n}}\text{=}\underset{n\rightarrow \infty }{\lim }\frac{\frac{1}{%
\beta -\alpha }\left[ \left( -b\text{+}a\beta \right) \alpha ^{n+1}\text{+}%
\left( b-a\alpha \right) \beta ^{n+1}\right] }{\frac{1}{\beta -\alpha }\left[
\left( -b\text{+}a\beta \right) \alpha ^{n}\text{+}\left( b-a\alpha \right)
\beta ^{n}\right] }\text{=}\max \{\alpha ,\beta \}.
\end{equation*}%
$\Box \medskip $

\textbf{Remark 3.7.}

1) Similar results can be obtained for the right Fibonacci type sequences.

2) The above formula $\left( 3.11\right) $ is a Binet's-type formula.
Indeed, for $A=B=1,a=0,b=1$, we obtain the Binet's formula for Fibonacci
sequence.

3) The limit $\mathcal{L}$, from relation $\left( 3.12\right) $ is a
Golden-ratio type number. Indeed, for $A=B=1,a=0,b=1$, we obtain the
Golden-ratio.\bigskip

We consider now the following difference equation of degree two%
\begin{equation}
d_{n}=ad_{n-1}+bd_{n-2},d_{0}=\alpha ,d_{1}=\beta ,  \tag{3.13.}
\end{equation}%
with $a,b,\alpha ,\beta $ arbitrary integers, $n\in \mathbb{N}$. It is clear
that $d_{2}=a\beta +b\alpha ,d_{3}=a^{2}\beta ++ab\alpha +b\beta ,$ etc.

Let $\left( G,\ast \right) $ be an arbitrary group and $g_{0},g_{1}\in G$.
We define the following sequence 
\begin{equation}
\varphi _{n}=\varphi _{n-1}^{a}\ast \varphi _{n-2}^{b},\varphi
_{0}=g_{1}^{\alpha }\ast g_{0}^{d_{-1}},\varphi _{1}=g_{1}^{\beta }\ast
g_{0}^{\alpha },  \tag{3.14.}
\end{equation}%
where $d_{-1}$ is defined depending on the sequence $\left( d_{n}\right)
_{n\in \mathbb{N}}$ and $g^{a}=g\ast g\ast ...\ast g,a-$time, $g\in G$. This
sequence is called the \textit{left} \textit{d-type sequence attached to the
binary relation }\textquotedblright $\ast $\textquotedblright , \textit{\
generated by }$g_{0}$ \textit{and} $g_{1}$.

In a similar way, we define the \textit{right} \textit{d-type sequence
attached to the binary relation }\textquotedblright $\ast $%
\textquotedblright , \textit{\ generated by }$g_{0}$ \textit{and} $g_{1}$,
namely 
\begin{equation}
\varphi _{n}=\varphi _{n-2}^{b}\ast \varphi _{n-1}^{a},\varphi
_{0}=g_{0}^{d_{-1}}\ast g_{1}^{\alpha },\varphi _{1}=g_{0}^{\alpha }\ast
g_{1}^{\beta },  \tag{3.15.}
\end{equation}%
where $d_{-1}$ is defined depending on the sequence $\left( d_{n}\right)
_{n\in \mathbb{N}}$ and $g^{a}=g\ast g\ast ...\ast g,a-$time, $g\in G$%
.\medskip 

\textbf{Proposition 3.8.} \textit{With the above notations, if }$g_{1}\ast
g_{0}=g_{0}\ast g_{1}$, \textit{the following relation holds:} 
\begin{equation}
\varphi _{n}=g_{1}^{d_{n}}\ast g_{0}^{d_{n-1}}.  \tag{3.16.}
\end{equation}

\textbf{Proof. }By straightforward calculations, it results\newline
$\varphi _{2}=\varphi _{1}^{a}\ast \varphi _{0}^{b}=\left( g_{1}^{\beta
}\ast g_{0}^{\alpha }\right) ^{a}\ast \left( g_{1}^{\alpha }\ast
g_{0}^{d_{-1}}\right) ^{b}=$\newline
$=g_{1}^{a\beta +b\alpha }\ast g_{0}^{a\alpha
+bd_{-1}}=g_{1}^{ad_{1}+bd_{0}}\ast g_{0}^{ad_{0}+bd_{-1}}=$\newline
$=g_{1}^{d_{2}}\ast g_{0}^{d_{1}}$\newline
$\varphi _{3}=\varphi _{2}^{a}\ast \varphi _{1}^{b}=\left( g_{1}^{d_{2}}\ast
g_{0}^{d_{1}}\right) ^{a}\ast \left( g_{1}^{d_{1}}\ast g_{0}^{d_{0}}\right)
^{b}=$\newline
$=g_{1}^{d_{3}}\ast g_{0}^{d_{2}}$. Using induction, we get\newline
$\varphi _{n+1}=\varphi _{n}^{a}\ast \varphi _{n-1}^{b}=\left(
g_{1}^{d_{n}}\ast g_{0}^{d_{n-1}}\right) ^{a}\ast \left( g_{1}^{d_{n-1}}\ast
g_{0}^{d_{n-2}}\right) ^{b}=$\newline
$=g_{1}^{d_{n+1}}\ast g_{0}^{d_{n}}.\Box \medskip $

\textbf{Remark 3.9. }

1) For $a=b=1,\alpha =0,\beta =1$, we get the Fibonacci sequence and
relation $\left( 3.14\right) $ becomes%
\begin{equation*}
\varphi _{n}=\varphi _{n-2}\ast \varphi _{n-1},\varphi
_{0}=g_{0}^{-1},\varphi _{1}=g_{1},
\end{equation*}%
since $f_{-1}=1$. In this case, relation $\left( 3.16\right) $ has the
following form%
\begin{equation}
\varphi _{n}=g_{1}^{f_{n}}\ast g_{0}^{f_{n-1}}.  \tag{3.17.}
\end{equation}%
\newline
Relation $\left( 3.17\right) $ was independently found in [HKN; 12],
Proposition 5.2 and it is a particular case of relation $\left( 3.16\right) $%
.

2) For $a=b=1,\alpha =2,\beta =1$, we get the Lucas sequence and relation $%
\left( 3.14\right) $ becomes 
\begin{equation*}
\varphi _{n}=\varphi _{n-2}\ast \varphi _{n-1},\varphi _{0}=g_{1}^{2}\ast
g_{0}^{-1},\varphi _{1}=g_{1}\ast g_{0}^{2}.
\end{equation*}%
In this case, relation $\left( 3.16\right) $ has the following form%
\begin{equation*}
\varphi _{n}=g_{1}^{l_{n}}\ast g_{0}^{l_{n-1}}.
\end{equation*}%
\begin{equation*}
\end{equation*}

\textbf{Conclusions.} In this paper, we provided new applications of
Fibonacci and Lucas sequences. In Proposition 3.2 and in Proposition 3.3, we
obtained some algebraic structures defined by using Fibonacci and Lucas
elements. In Proposition 3.5, we provided some new properties and relations
of these sequences.

In Proposition 3.6 and in Proposition 3.8, we also generalized Fibonacci and
Lucas numbers by using an arbitrary binary relation over the real fields
instead of addition of the real numbers and we give properties of the new
obtained sequences.

Moreover, taking into consideration some known relations between these
numbers, in Proposition 3.4, we give a method to find new examples of
quaternion split algebras. This method can be an easy alternative for using
the properties of quadratic forms.

As a further research, we intend to continue the study of such sequences,
especially in their connections with quaternion algebras, expecting to
obtain new and interesting results.

\begin{equation*}
\end{equation*}

\textbf{References} \bigskip \newline
\textbf{[}Fib.\textbf{]}
http://www.maths.surrey.ac.uk/hosted-sites/R.Knott/Fibonacci/fib.html\newline
\textbf{[}Fl, Sa; 15\textbf{]} C. Flaut, D. Savin, \textit{Quaternion
Algebras and Generalized Fibonacci-Lucas Quaternions}, Adv. Appl. Clifford
Algebras, 25(4)(2015), 853-862.\newline
\textbf{[}Fl, Sh; 13\textbf{]} C. Flaut, V. Shpakivskyi, \textit{On
Generalized Fibonacci Quaternions and Fibonacci-Narayana Quaternions}, Adv.
Appl. Clifford Algebras, 23(3)(2013), 673--688.\newline
\textbf{[}Gi, Sz; 06\textbf{]} Gille, P., Szamuely, T., Central Simple
Algebras and Galois Cohomology, Cambridge University Press, 2006.\newline
[HKN; 12] Han, J.S., Kim, H.S., Neggers, J., \textit{Fibonacci sequences in
groupoids}, Advances in Difference Equations 2012, 2012:19, \newline
\textbf{[}Ho; 63\textbf{]}\ A. F. Horadam, \textit{Complex Fibonacci Numbers
and Fibonacci Quaternions}, Amer. Math. Monthly, 70(1963), 289--291.\newline
\textbf{[}Mag.\textbf{]} http://magma.maths.usyd.edu.au/magma/handbook/ 
\newline
\textbf{[}Sa; 19\textbf{]} D. Savin, Special numbers, special quaternions
and special symbol elements, chapter in the book Models and Theories in
Social Systems, vol. 179, Springer 2019, ISBN-978-3-030-00083-7 , p.
417-430. 
\begin{equation*}
\end{equation*}

Cristina FLAUT

{\small Faculty of Mathematics and Computer Science, }\newline

{\small Ovidius University of Constan\c{t}a, Rom\^{a}nia,}

{\small Bd. Mamaia 124, 900527,}

{\small http://www.univ-ovidius.ro/math/}

{\small e-mail: cflaut@univ-ovidius.ro; cristina\_flaut@yahoo.com}%
\begin{equation*}
\end{equation*}%
\medskip \qquad\ \qquad\ \ 

Diana SAVIN

{\small Faculty of Mathematics and Computer Science, }

{\small Ovidius University of Constan\c{t}a, Rom\^{a}nia, }

{\small Bd. Mamaia 124, 900527, }

{\small http://www.univ-ovidius.ro/math/}

{\small e-mail: \ savin.diana@univ-ovidius.ro, \ dianet72@yahoo.com}

\begin{equation*}
\end{equation*}%
\qquad

Geanina ZAHARIA

{\small PhD student at Doctoral School of Mathematics,}

{\small Ovidius University of Constan\c{t}a, Rom\^{a}nia,}

{\small geaninazaharia@yahoo.com}

\end{document}